\def\IR{{\Bbb R}} 
 
\def\IS{{\Bbb S}} 
 
\def\IK{{\Bbb K}}

\def\IC{\Bbb C} 
\def\ID{{\Bbb D}}

\def\zbar{{\overline{z}}} 
 
\def\wbar{{\overline{w}}}

\documentclass[10pt]{article} 

\usepackage[psamsfonts]{amssymb} 
\usepackage{amsfonts,amsmath}  

\usepackage{epsfig,multicol}

\newtheorem{theorem}{Theorem}[section] 
 
\newtheorem{lemma}{Lemma}[section] 
\newtheorem{corollary}{Corollary}[section]

\newtheorem{conjecture}{Conjecture}[section]

\numberwithin{equation}{section}

\title{Higher regularity and uniqueness for inner variational equations.}
\author{Gaven Martin \& Cong Yao \thanks{ This research of both authors is supported in part by a grant from the NZ Marsden Fund.\newline
\newline
{\bf  Mathematics subject classification 2010,}{Primary: 30C62 31A05 49J10;  } \newline 
{\bf  Keywords,} {Calculus of variations,  quasiconformal, distributional equations, mean distortion}}}
\date{}
\begin{document}
\maketitle

\begin{abstract}  
We study local minima of the $p$-conformal energy functionals,
\[
\mathsf{E}_{\cal A}^\ast(h):=\int_\ID {\cal A}(\IK(w,h)) \;J(w,h) \; dw,\quad h|_\IS=h_0|_\IS,
\]
defined for self mappings $h:\ID\to\ID$ with finite distortion of the unit disk with prescribed boundary values $h_0$.
Here $\IK(w,h) = \frac{\|Dh(w)\|^2}{J(w,h)} $ 
is the pointwise distortion functional,  and ${\cal A}:[1,\infty)\to [1,\infty)$ is convex and increasing with ${\cal A}(t)\approx t^p$ for some $p\geq 1$,  with additional minor technical conditions. Note ${\cal A}(t)=t$ is the Dirichlet energy functional.

Critical points of $\mathsf{E}_{\cal A}^\ast$  satisfy the Ahlfors-Hopf inner-variational equation 
\[ {\cal A}'(\IK(w,h)) h_w \overline{h_\wbar} = \Phi \]
where $\Phi$ is a holomorphic function.  Iwaniec, Kovalev and Onninen established the Lipschitz regularity of critical points. Here we give a sufficient condition to ensure that a local minimum is a diffeomorphic  solution to this equation,  and that it is unique.  This condition is necessarily satisfied by any locally quasiconformal critical point,  and is basically the assumption  $\IK(w,h)\in L^1(\ID)\cap L^r_{loc}(\ID)$ for some $r>1$.  
\end{abstract}

\section{Introduction}

A mapping $f:\ID\to \ID$ has finite distortion if 
\begin{enumerate}
\item $f\in W^{1,1}_{loc}(\ID)$,  the Sobolev space of functions with locally integrable first derivatives,
\item the Jacobian determinant $J(z,f)\in L^{1}_{loc}(\ID)$, and 
\item there is a measurable function ${\bf K}(z)\geq 1$, finite almost everywhere, such that 
 \begin{equation}\label{1.1}
 |Df(z)|^2 \leq {\bf K}(z) \, J(z,f), \hskip10pt \mbox{ almost everywhere in $\ID$}.
 \end{equation}
\end{enumerate}
See \cite[Chapter 20]{AIM} or \cite{HK2} for the basic theory of mappings of finite distortion and the associated governing equations; degenerate elliptic Beltrami systems.
In (\ref{1.1}) the operator norm is used.  However this norm loses smoothness at crossings of the singular-values of the differential $Df$ and for this reason when considering minimisers of distortion functionals one considers the distortion functional
\begin{equation}\label{1.2}
\IK(z,f) = \left\{\begin{array}{cc} \frac{\|Df(z)\|^2}{J(z,f)}, & \mbox{if $J(z,f)\neq 0$} \\
1, & \mbox{if $J(z,f)= 0$.} \end{array}\right.
\end{equation}
This was already realised by Ahlfors in his seminal work proving Teichm\"uller's theorem and establishing the basics of the theory of quasiconformal mappings, \cite[\S 3, pg 44]{Ahlfors}. We reconcile (\ref{1.1}) and (\ref{1.2}) by noting $\IK(z,f)=\frac{1}{2}\big({\bf K}(z)+1/{\bf K}(z) \big)$ almost everywhere,  where ${\bf K}(z)$ is chosen to be the smallest functions such that (\ref{1.1}) holds.
 
 \medskip
 
\noindent Let ${\cal A}:[1,\infty)\to [1,\infty)$ be convex and increasing with 
\begin{equation}\label{Acondition}
 p{\cal A}(t) \leq t {\cal A}'(t), \hskip15pt \mbox{ for some $p>1.$}
\end{equation}
The number $p$ here determines the higher regularity assumptions we make.
 The ${\cal A}$-mean distortion of a self-homeomorphism of $\overline{\ID}$  is defined as
\begin{equation}\label{1.3}
\mathsf{E}_{\cal A}(f):=\int_\ID {\cal A}(\IK(z,f))\; dz,
\end{equation}
The canonical examples are when ${\cal A}(t)=t^p$ and there we simply write $\mathsf{E}_{p}(f)$.  The dual energy functional is 
\begin{equation}\label{1.4}
\mathsf{E}_{\cal A}^\ast(h):=\int_\ID {\cal A}(\IK(w,h))\; J(w,h) \; dw,
\end{equation}
For self homeomorphisms of $\ID$,  $f$ and $h=f^{-1}$, of finite distortion we have
\begin{equation}
\mathsf{E}_{\cal A}^\ast(h) = \mathsf{E}_{\cal A}(f).
\end{equation}
See \cite{H} or \cite{AIM,HK1,HK2} for more information on the change of variables needed here.\\

We recall the following conjecture in \cite{IMO}.
\begin{conjecture}
Let $f_0:\overline{\ID}\to\overline{\ID}$ be a homeomorphism of finite distortion with $\mathsf{E}_{\cal A}(f_0)<\infty$.  In the space of homeomorphic mappings of finite distortion with boundary values $f_0$, there is a minimiser $f$ which is also a smooth diffeomorphism.
\end{conjecture} 

There is of course a similar conjecture for $h$ and either one implies the other.

\subsection{Inner variational equations.}

Note that the a priori regularity for $f$ in (\ref{1.3}) is $W^{1,\frac{2p}{p+1}}(\ID)$ and for $h$ in (\ref{1.4}) is $W^{1,2}(\ID)$.   Let $\varphi\in C^{\infty}_{0}(\ID)$ with $\|\nabla\varphi\|_{L^{\infty}(\ID)}<\frac{1}{2}$.  Then for $t\in (-1,1)$ the mapping $g^t(z)=z+t\varphi(z)$ is a diffeomorphism of $\ID$ to itself which extends to the identify on the boundary $\IS$.  If $f$ is a mapping of finite distortion for which $\mathsf{E}_{\cal A}(f)<\infty$,  then so is $\mathsf{E}_{\cal A}(f\circ g^t)<\infty$ and they share boundary values. Similarly for $h$ and $\mathsf{E}_{\cal A}^\ast(h\circ g^t)$.

The functions $t\mapsto \mathsf{E}_{\cal A}(f\circ g^t)$ and $t\mapsto \mathsf{E}_{\cal A}^\ast(h\circ g^t)$ are smooth function of $t$. Thus if  $f$ or $h$ is a minimiser in any reasonable class we have the stationary equations
\begin{eqnarray*}
\frac{d}{dt} \Big|_{t=0} \mathsf{E}_{\cal A}(f\circ g^t)=0, &&
\frac{d}{dt} \Big|_{t=0} \mathsf{E}_{\cal A}^\ast(h\circ g^t) = 0.
 \end{eqnarray*}
It is a calculation to verify that the first equation is equivalent to 
\begin{equation}\label{1.5}
2p\int_\ID \IK_f {\cal A}'(\IK) \frac{\overline{\mu_f}}{1+|\mu_f|^2}\varphi_\zbar dz=\int_\ID {\cal A}(\IK) \varphi_zdz,\quad\forall\varphi\in C_0^\infty(\ID ).
\end{equation}
and that the second is equivalent to 
\begin{equation}\label{1.6}
{\cal A}'(\IK(w,h))h_w \overline{h_\wbar} = \Phi.
\end{equation}
where $\Phi$ is holomorphic.  Using an early modulus of continuity estimate  Alhfors showed that there is always an $h:\overline{\ID}\to\overline{\ID}$ with $h\in C(\overline{\ID})\cap W^{1,2}(\ID)$ minimising (\ref{1.4}) for quasisymmetric boundary values and therefore solving (\ref{1.6}).  Strictly speaking he used ${\cal A}(t)=t^p$, $p\geq 2$, but more recent modulus of continuity estimates give the more general result,  \cite{AIM,GV,HK2}.   For this reason we call $\Phi$ the Ahlfors-Hopf differential.\\

The two strongest results currently known to us are the Lipschitz regularity of Iwaniec, Kovalev and Onninen \cite{IKO}
\begin{theorem} \label{loclip}
Let $h\in W^{1,2}(\ID)$ be a mapping of finite distortion which solves (\ref{1.6}) for holomorphic $\Phi$.  Then $h$ is locally Lipschitz.
\end{theorem}
Also our earlier result (strictly speaking only for ${\cal A}(t)=t^p$ but the ideas are exactly the same) \cite{MY}.
\begin{theorem}\label{thm1.1}
Let $f$ be a finite distortion function that satisfies the distributional equation  (\ref{1.5}). Assume that
$\mathbb{K}(z,f) \in L_{loc}^r(\ID)$,  for some $r>p+1$. Then $f$ is a local diffeomorphism in $\ID$.
\end{theorem}
Actually $r$ does not have to be uniform in $\ID$.  The following corollary is almost immediate.

\begin{corollary}
Let $h:\overline{\ID} \to\overline{\ID }$  be a continuous locally quasiconformal solution to  (\ref{1.6}). Then $h$ is a smooth self-diffeomorphism of $ \ID$.
\end{corollary}

We also gave a counterexample to justify some assumptions on the integrability of the distortion.
\begin{theorem}\label{thm 1.2}
There is a Sobolev mapping with $f\in W^{1,\frac{2p}{p+1}}(\ID)$ with Beltrami coefficient $\mu_f$ satisfying the distributional equation (\ref{1.5}),  and with $\IK_f \in L^p(\ID)\setminus\bigcup_{q>p} L^q_{loc}(\ID)$.
In particular,  this mapping $f$ has $\mathsf{E}_p(f)<\infty$ and solves the distributional equation, but it cannot be locally quasiconformal.
\end{theorem}
The mapping $f$ of Theorem \ref{thm 1.2} has a pseudo-inverse $h:\overline{\ID}\to\overline{\ID}$,  $h\in C(\ID)\cap W^{1,2}_{loc}(\ID)$, a monotone mapping for which $h(f(z))=z$ for almost every $z\in \ID$.

Unfortunately we do not know if $h$ is a homeomorphism, if $h$ has homeomorphic boundary values or even if $h$ satisfies the Ahlfors-Hopf equation (though this last would follow if $f$ were a local minimum for (\ref{1.4})).  In the case $p=1$ and ${\cal A}(t)=t$, \cite[Example 3.4]{IKO} provides a Lipschitz solution to the Hopf-Laplace equation $h_w\overline{h_\wbar}=-1$ with $J(w,h)\geq 0$ almost everywhere and yet is not a homeomorphism.  This map can be modified so as to be defined on $\ID$,  but its image seems unwilling to be modified so as to be a disk without avoiding the singular set and so becoming a diffeomorphism. 

We are unaware of a way to connect the two inner-variational equations,  even for homeomorphic solutions,  without fairly strong a priori assumptions.  As for the boundary values,  it remains unclear as to what is the exact criterion for a self-homeomorphism $f_0$ or $h_0$ to admit an appropriate extension of finite energy,  so that the family we might consider is not empty (though see \cite{AIMO} in the case $p=1$).  Thus we stick to the class of quasisymmetric mappings which have a quasiconformal extension.  Our main results here are the following.{\em We shall always assume that $p>1$ unless otherwise stated}.

\begin{theorem}\label{main1} Let $h:\ID\to\ID$ with $\mathsf{E}_{\cal A}^\ast(h)<\infty$ for quasisymmetric boundary values $h_0:\IS\to\IS$.  If $\IK(w,h)\in L^1(\ID)$,  then $h$ is a homeomorphism. If in addition $h$ is a local minimum for $\mathsf{E}_{\cal A}^\ast$ and if $\IK(w,h)\in L^r_{loc}(\ID)$ for some $r>1$,  then $h$ is a diffeomorphism.
\end{theorem}

We remark that Iwaniec's calculation of the second inner-variation (personal communication) suggests that there are in fact no local maxima.  We note a slightly different result.  That $h$ is a homeomorphic local maximum or minimum shows its inverse satisfies the inner distributional equation.  This would also be guaranteed if $h$ satisfies the outer distributional equation
\[ \frac{d}{dt} \int_\ID \IK^p(w,g^t\circ h)\, J(w,g^t\circ h) \,dw = 0\]
It is a lengthy calculation to reveal this equation is
\begin{equation}  
 \int_\ID \IK^p_h \big((\IK_h+1)p-1 \big)h_\wbar \phi_w \; dw   = \int_\ID \IK^p_h \big((\IK_h-1)p-1 \big)h_w \phi_\wbar \; dw \label{1.8}
\end{equation}
A finite distortion mapping which is a solution to (\ref{1.8})  and has $\mathsf{E}_{\cal A}^\ast(h)<\infty$ is called an {\em outer variational stationary point}.
\begin{theorem}\label{main2}  If $h$ is an outer-variational stationary point with quasisymmetric boundary values,  if  $\IK(w,h)\in L^1(\ID)$  and if $E_{q}^\ast(h)$ is locally finite for some $q>p+1$,  then $h$ is a diffeomorphism.
\end{theorem}
This is in essence a restatement of Theorem \ref{thm1.1}.  \\

We can say a little about uniqueness here too.

\begin{theorem}\label{unique} Let $h:\overline{\ID}\to\overline{\ID}$, $h\in W^{1,2}(\ID)$, be a diffeomorphism of $\ID$ with homeomorphic boundary values, and with Ahlfors-Hopf differential  $\Phi$ as at (\ref{1.6}). 
Let $g:\overline{\ID}\to\overline{\ID}$ be continuous and a homeomorphism on $\IS$.  Suppose $g$ is a solution to the Ahlfors-Hopf equation 
\begin{equation} 
{\cal A}'(\IK(w,g))g_w \overline{g_\wbar} = \Phi,
\end{equation}
in $\ID$ with $g(0)=h(0)$ and $g(1)=h(1)$.  Then $g\equiv h$.
\end{theorem}

Note that we do not require that $g=h$ on $\IS$ in the hypotheses. 
In fact diffeomorphic minimisers are locally absolute minimisers for their boundary values.

\begin{theorem} Let $\Omega_1,\Omega_2$ be Jordan domains and let $h:\Omega_1\to\Omega_2$ be a diffeomorphic minimiser of $\mathsf{E}_{\cal A}^\ast$ for its boundary values.  If   $D = D(z_0,r)\subset\Omega_2$ is any disk and  $\varphi:h^{-1}(D)\to \ID$ is a Riemann map, then $h_\ast(z) = \frac{1}{r} (h\circ \phi^{-1})(z)-z_0 :\ID\to \ID$ is the unique minimiser for its boundary values.
\end{theorem}
\noindent{\bf Proof.}  Suppose there is $g$ with the same or smaller energy in $\ID$ for the boundary values of $h_\ast$.  Set
\begin{equation}
\tilde{h} = \left\{ \begin{array}{ll} h(z), & z\in \ID\setminus \{h^{-1}(D)\}, \\
 rg(\varphi(z))+z_0, & z\in  \{h^{-1}(D)\}  \end{array} \right. .
\end{equation}
Then $\tilde{h}$ is a mapping of finite distortion  and has energy no more than $h$,  a minimiser (the only issue is the $W^{1,1}_{loc}$ regularity across the set $\partial h^{-1}(D)$,  but this is a smooth Jordan curve).  So $\tilde{h}$ is also a minimiser,  and therefore has a holomorphic Ahlfors-Hopf differential.  This differential must be $\Phi$ since $\tilde{h}$ it agrees with $h$ on $\ID\setminus h^{-1}(D)$.  The result now follows from Theorem \ref{unique}. \hfill $\Box$

\section{Diffeomorphisms; proof of Theorem \ref{main1}.}
 
Since ${\cal A}(t) \geq t$ we see $h$  lies in the Sobolev space $W^{1,2}(\ID)$ as  $\mathsf{E}_{\cal A}^\ast(h)<\infty$.  We observe that since $p>1$, the estimate
\[  \frac{\|Dh\|^{2p-2}}{J(w,h)^{p-1}} \frac{|\mu|}{1+|\mu|^2} \approx |\Phi| \]
implies that $h$ is a mapping of finite distortion.

Let $h_\ast:\IC\setminus \ID \to \IC\setminus \ID$ be a quasiconformal extension of the quasisymmetric boundary values $h_0$. We can ensure that $h_\ast(z)=z$ for all sufficiently large $z$ by the quasiconformal version of the Sch\"onflies theorem,  \cite{GVam} or \cite[\S 7]{GMP}. Then 
\[ H(w) = \left\{\begin{array}{ll} h(w), & w\in \ID, \\ h_\ast(w), &w\in \IC\setminus \ID, \end{array} \right. \]
is a mapping of finite distortion,  $H(w)-w\in W^{1,2}(\IC)$ and $\mu_H=H_\wbar/H_w$ is compactly supported.  Further $\IK(w,H)-1\in L^1(\IC)$.  Then \cite[Theorem 20.2.1]{AIM} provides a homeomorphic entire principal solution (we have to make the minor adjustment of replacing $\ID$ by a larger disk in which $\mu$ is compactly supported) to the Beltrami equation $g_\wbar=\mu_Hg_w$ and the Stoilow factorisation shows this solution must be $H$ up to a similarity.  Thus $H$,  and hence $h$, is a homeomorphism.  

Next we suppose that $h$ is a homeomorphic local maximum or local minimum and set $f=h^{-1}:\ID\to\ID$.  We have $\mathsf{E}_{\cal A}^\ast(h)=\mathsf{E}_{\cal A}(f)$ and also,  for $g^t(z)=z+t\varphi(z)$ with compactly supported test function $\varphi$, $|\nabla\varphi|<1$, 
\begin{equation}
\mathsf{E}_{\cal A}(f\circ g^t)= \mathsf{E}_{\cal A}^\ast((g^t)^{-1}h).
\end{equation}
Then $h$ a local  minimum for $\mathsf{E}_{\cal A}^\ast$ implies that $\frac{d}{dt}\big|_{t=0} \mathsf{E}_{\cal A}(f\circ g^t) = 0$ and $f$ satisfies the distributional equation (\ref{1.5}).  Next,  the Alhfors-Hopf differential is
\begin{eqnarray*} 
\Phi & = & {\cal A}'(\IK(w,h)) h_w \overline{h_\wbar} = {\cal A}'(\IK(w,h)) |h_w|^2 \overline{\mu_h} \\ & = & \IK(w,h) {\cal A}'(\IK(w,h)) J(w,h) \frac{\overline{\mu_h}}{1+|\mu_h|^2} \\
|\Phi |&\geq  & c_0\, \IK(w,h)^p  J(w,h) \end{eqnarray*}
The last inequality holding by virtue of (\ref{Acondition}) and only on $E=\{w:|\mu_h(w)|\geq \frac{1}{2}\}$ for $c_0$ a positive constant. In particular we see that  $\IK(w,h)^p  J(w,h)$ is locally bounded on the set $\IK(w,h)\geq \frac{5}{3}$.  If $\IK(w,h)\in L^{r}(V)$ for some relatively compact set $V$,  then  $\IK(w,h)^{p+r}  J(w,h) \in L^1(V)$ as the Jacobian is locally integrable.  Then
\[ \int_{h^{-1}(V) }\IK(z,f)^{p+r}(z,f)  = \int_\ID \IK(w,h)^{p+r}  J(w,h) <\infty \]
and since $h:\ID\to\ID$ is a homeomorphism we see our hypotheses imply $\IK(z,f)\in L^{q}_{loc}(\ID)$ for some $q>1+p$ and hence $f:\ID\to\ID$ is a diffeomorphism by Theorem \ref{thm1.1}. The result now follows. \hfill $\Box$
   
\section{Uniqueness: Proof of Theorem \ref{unique}}
 
The proof of uniqueness will be separated into two parts.  First we will find   a degenerate elliptic Beltrami equation so that $h$ has the Hopf differential $\Phi$ if and only if 
\begin{equation}\label{3.1}
h_\wbar=B(w,h_w).
\end{equation}
Here $B(w,\xi):\ID\times\IC\to\IC$ is defined implicitly and is smooth away from the set $\xi=0$.  We then give the ellipticity bounds on the nonlinear equation (\ref{3.1}).  We discuss the Schauder bounds and smoothness elsewhere.

Second, we use the ellipticity bounds, together with the existence of a diffeomorphic solution to establish the following lemma.   

\begin{lemma}\label{lem1}
Let $B(w,\xi)$ as above and $h:\overline{\ID}\to\overline{\ID}$ a continuous $W^{1,2}(\ID)$ solution to (\ref{3.1}) and homeomorphic on $\IS$.  Let $g$ be a homeomorphism from $\overline{\ID}$ to $\overline{\ID}$, a diffeomorphism from $\ID$ to $\ID$, lies in $W^{1,2}(\ID)$ and that also satisfies equation (\ref{3.1}). Then $\eta:=g-h$ is a locally quasiregular mapping.
\end{lemma}

Given this lemma, uniqueness quickly follows in exactly the same way as in \cite[\S 9.2.2, pp 267]{AIM} using the total variation of the boundary values \cite[Lemma 9.2.2]{AIM} and the the Stoilow factorisation theorem,  as in \cite[Lemma 9.2.3]{AIM} .

\subsection{The equation;  $B(w,\xi)$ and its ellipticity properties}
We may assume that $\Phi$ is not identically $0$. We make the simplifying assumption that ${\cal A}(t)=t^p$.  The ideas in the general case are the same,  we consider the level curves of the function $(x,y)\mapsto {\cal A}'((x^2+y^2)/(x^2-y^2)) \, x y$,  but some of the formulas become a little unwieldy and estimates not as clean.

\medskip

We begin by considering the level curves of the function
\[
(x,y)\mapsto \Big(\frac{x^2+y^2}{x^2-y^2}\Big)^{p-1} \, x y.
\]
This function is defined on the region $\Omega=\{(x,y)\in \IR^2: x>0, 0<y<x\}$.
For fixed $x>0$ the function $y\mapsto(\frac{x^2+y^2}{x^2-y^2})^{p-1} \, x y $ is strictly increasing for $0<y<x$:
\[
\frac{\partial}{\partial y}\Big[\Big(\frac{x^2+y^2}{x^2-y^2}\Big)^{p-1}xy\Big]=\Big(\frac{x^2+y^2}{x^2-y^2}\Big)^{p-1}x+(p-1)\Big(\frac{x^2+y^2}{x^2-y^2}\Big)^{p-2}\frac{4x^3y^2}{(x^2-y^2)^2}>0.
\]
Hence for any $k>0$ there is a unique solution $y$ so that 
\[
\Big(\frac{x^2+y^2}{x^2-y^2}\Big)^{p-1}x y = k.
\]
The implicit function theorem guarantees the level curves are simple arcs in the $x,y$ plane and on the level curve we have  that $y$ can be expressed as a function of $x$,
$ y = A_k(x)$.

{ 
\scalebox{0.4}{\includegraphics[viewport= -50 480 620 780]{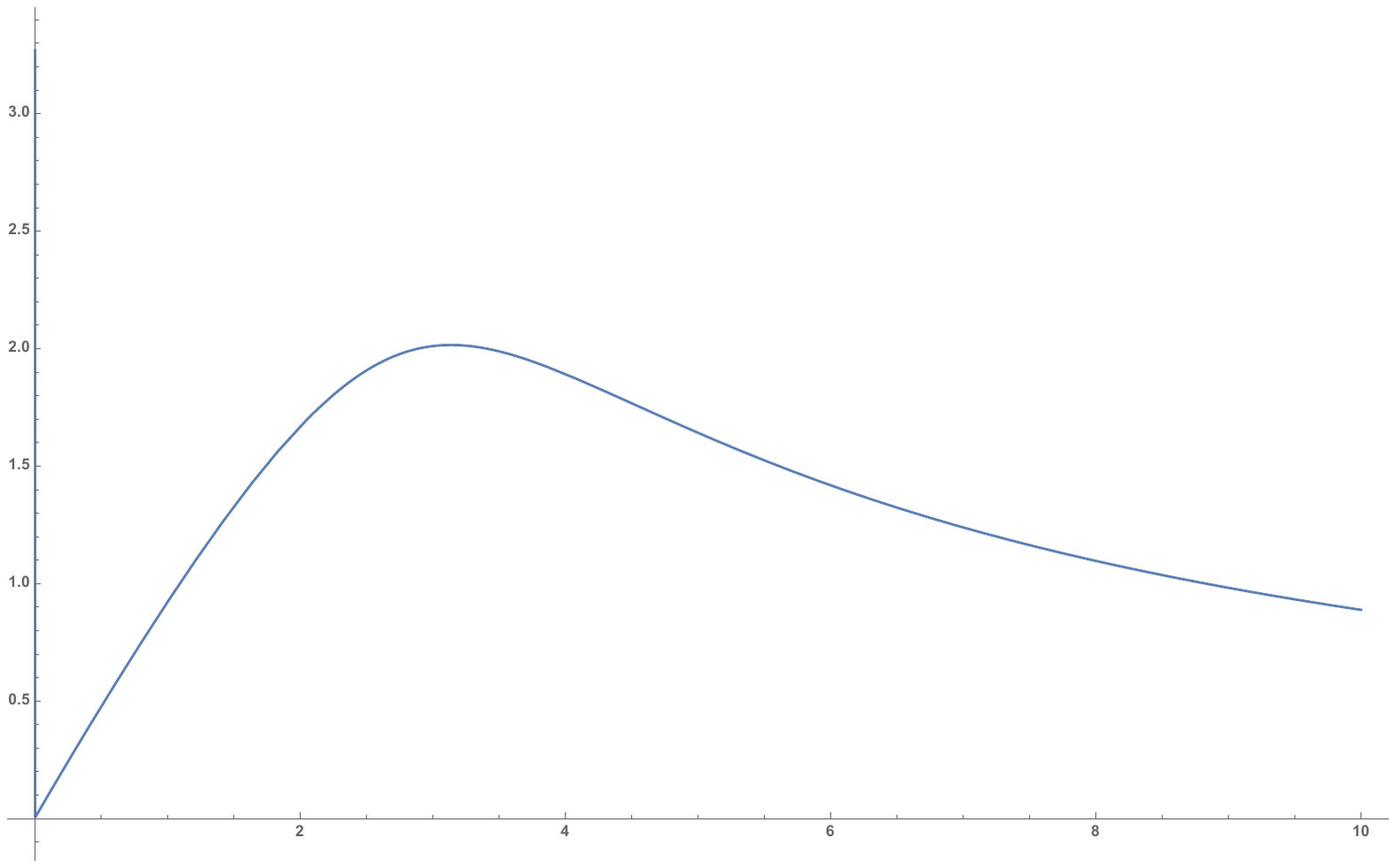}}
} 

\noindent {\bf Figure 1.} The graph of the level curve $\frac{x^2+y^2}{x^2-y^2}\, x y = 10$.

\medskip
 
On this curve  
\begin{equation}\label{3.2}
(p-1)\log[\frac{x^2+A_k^2(x)}{x^2-A_k^2(x)}]+\log x+\log A_k(x)=\log k.
\end{equation}
With $x=|h_w|$ and $y=|h_\wbar|$ we have
\begin{eqnarray*}
|h_\wbar| = A_k(|h_w|), &&
h_\wbar\overline{h_w}\frac{\Phi}{|\Phi|} = A_k(|h_w|)|h_w|,
\end{eqnarray*}
and hence we find a nonlinear Beltrami equation for $h$ as
\begin{equation}\label{3.3}
h_\wbar   =\frac{\overline{\Phi}}{|\Phi|}  A_{|\Phi(w)|}(|h_w|) \frac{h_w}{|h_w|} = B(w,h_w).
\end{equation}
Next,  the ellipticity properties of equation (\ref{3.3}). We drop  the subscript on $A$, 
\begin{eqnarray}\label{3.4}
|B(w,\zeta)-B(w,\xi)| & = &  \left|A(|\zeta|) \frac{\zeta}{|\zeta|} -A(|\xi|) \frac{\xi}{|\xi|} \right|.
\end{eqnarray}
We set $V(x)=A(x)/x$. After division, equation (\ref{3.4}) reads as 
\begin{eqnarray}
\frac{|B(w,\zeta)-B(w,\xi)| }{|\zeta-\xi|}& = & \frac{ \big|V(|\zeta|)  \zeta -V(|\xi|)  \xi \big|}{|\zeta-\xi|}.
\end{eqnarray}
 We put  $|\zeta|=t$,  $|\xi|=s$, $a=V(t)$ and $b=V(s)$.  Then there is a $\theta\in[0,2\pi]$ such that $
\zeta\cdot\xi= st\cos(\theta)$, 
and
\begin{eqnarray}\label{3.6}
\frac{ \big|V(|\zeta|)  \zeta -V(|\xi|)  \xi \big|^2}{|\zeta-\xi|^2} & = & \frac{a^2 t^2+b^2 s^2 - 2 abst \cos(\theta)}{t^2+s^2-2st\cos(\theta)}:=F(\theta).
\end{eqnarray}
 We differentiate (\ref{3.6}) with respect to $\theta$ to see
\[
\frac{d}{d\theta}F(\theta)=\frac{2  st[abt^2+abs^2-a^2 t^2-b^2 s^2 ]}{(t^2+s^2-2st\cos(\theta))^2} \sin(\theta) = \frac{2  st(a-b)(s^2b-t^2a)}{(t^2+s^2-2st\cos(\theta))^2}\sin(\theta).
\]
Here we claim that
\begin{equation}\label{3.7}
(a-b)(s^2b-t^2a)\geq0.
\end{equation}
Recall $V(x)=A(x)/x$. We also define $W(x)=xA(x)$. Then (\ref{3.2}) gives us the relation
\begin{eqnarray*}
\log k & = & (p-1)\log[\frac{1+V^2(x)}{1-V^2(x)}]+ 2 \log x +\log V(x), \\
  & = & (p-1)\log[\frac{x^4+W^2(x)}{x^4-W^2(x)}]+\log W(x),
\end{eqnarray*}
which we differentiate to see that
\[
V^\prime(x) \left[  \frac{4(p-1)V(x)}{1-V^4(x)}  + \frac{1}{V(x)}\right] = - \frac{2}{x},
\]
\[
W^\prime(x)\left[  \frac{4(p-1)x^4W(x)}{x^8-W^4(x)}  + \frac{1}{W(x)}\right] =\frac{8(p-1)x^3W^2(x)}{x^8-W^4(x)}.
\]
So $V$ is decreasing and $W$ is increasing. Now assume $t\leq s$, then
\[
a=V(t)\geq V(s)=b,
\hskip10pt
s^2b=W(s)\geq W(t)=t^2a,
\]
and vice versa. So (\ref{3.7}) follows. Assume $\zeta\neq\xi$, we then have
\[
\frac{d}{d\theta}\frac{\big|V(|\zeta|)\zeta -V(|\xi|)\xi\big|^2}{|\zeta-\xi|^2}=G(|\zeta|,|\xi|,\cos(\theta))\sin\theta,
\]
where $G(|\zeta|,|\xi|,\cos(\theta))$ is always non-negative. Then, in a period $\theta\in[0,2\pi]$, $F(\theta)$ is increasing when $\theta$ is moving from $0$ to $\pi$ and decreasing when $\theta$ is moving from $\pi$ to $2\pi$, so we get the maximum of $F(\theta)$ at $\theta=\pi$. In particular we can now write (\ref{3.4}) as
\begin{eqnarray}
\frac{|B(w,\zeta)-B(w,\xi)| }{|\zeta-\xi|}&\leq & \frac{ \big|V(t) t +V(s)s\big|}{|t+s|} = \frac{ \big|A(t) +A(s)\big|}{|t+s|}\nonumber  \\
& \leq & \max\{\frac{A(|\zeta|)}{|\zeta|},\frac{A(\xi)}{|\xi|}\}, \label{3.9}
\end{eqnarray}
whenever $\zeta\neq\xi$.\\

We now assume that $h$, $g$ are finite distortion homeomorphisms solutions to (\ref{3.3}) and consider the function $\eta=g-h$. 

At all but a discrete set of points $w\in\Omega$, we have $\Phi(w)\neq 0$. If $h_w=g_w$, then (\ref{3.3}) gives $h_\wbar=g_\wbar$; if $h_w\neq g_w$, then by (\ref{3.9}),
\[
|\mu_\eta|=\Big|\frac{\eta_\wbar}{\eta_w}\Big|=\Big|\frac{g_\wbar-h_\wbar}{h_w-g_w}\Big|\leq\max\{|\mu_g|,|\mu_h|\}.
\]
Also note
\begin{eqnarray*} |\eta_w|^2-|\eta_\wbar|^2 &=& |g_w-h_w|^2-|g_\wbar-h_\wbar|^2 \\
&=&  J(z,g) +J(z,h) -2 \Re e[g_w \overline{h_w} - g_\wbar \overline{h_\wbar}]\in L^1(\Omega).
\end{eqnarray*}
These facts imply that $\eta$ is a finite distortion function.\\

If we assume further that $g$ is diffeomorphic from $\ID$ to $\ID$, then in any compact subset $\Omega\subset\subset\ID$, we have
\[
|\Phi|\leq M<\infty,\quad|g_\wbar|\geq\varepsilon>0,\quad|\mu_g|\leq k<1.
\]
Now at any point $w\in\Omega$, by (\ref{3.9}),
\[
\frac{|h_\wbar-g_\wbar|}{|h_w-g_w|}\leq\frac{|h_\wbar|+|g_\wbar|}{|h_w|+|g_w|}\leq\frac{|h_\wbar|+|g_\wbar|}{|g_w|}.
\]
So we can choose a $\delta>0$ such that if $|h_\wbar|<\delta$, then
\[
\frac{|h_\wbar-g_\wbar|}{|h_w-g_w|}\leq\frac{|g_\wbar|}{|g_w|}+\frac{1-k}{2}<1.
\]
Note this $\delta$ depends only on $\varepsilon$, $M$ and $k$ but not a specific point $w\in\Omega$. On the other hand, if $|h_\wbar|>\delta$, then
\[
\IK_h^{p-1}|h_wh_\wbar|=|\Phi|\leq M,
\]
which gives 
$
\IK_h<(\frac{M}{\delta^2})^\frac{1}{p-1},
$ 
so at this point we have
\[
|\mu_\eta|\leq\max\Big\{\frac{1+|\mu_g|}{2},\sqrt{\frac{(\frac{M}{\delta^2})^\frac{1}{p-1}-1}{(\frac{M}{\delta^2})^\frac{1}{p-1}+1}}\Big\} < 1.
\]
This estimate holds locally uniformly and thus proves that $\eta$ is locally quasiregular and completes the proof of Lemma \ref{lem1}. \hfill $\Box$

\medskip

GJM Institute for Advanced Study, Massey University, Auckland, New Zealand \\
email: G.J.Martin@Massey.ac.nz

CY Institute for Advanced Study, Massey University, Auckland, New Zealand \\
email: C.Yao@massey.ac.nz 
\end{document}